\documentclass{amsart}
\usepackage[all]{xy}
\usepackage{graphicx}
\usepackage{indentfirst}
\usepackage{bm}
\usepackage{amsthm}
\usepackage{mathrsfs}
\usepackage{latexsym}
\usepackage{amsmath}
\usepackage{amssymb}
\usepackage{hyperref}
\usepackage{microtype}
\usepackage{color}
\usepackage{tikz}
\usepackage{comment}
\usepackage{appendix}
\usepackage{bookmark}
\allowdisplaybreaks

\theoremstyle{plain}
\newtheorem{theorem}{Theorem~}[section]

\newtheorem{lemma}[theorem]{Lemma~}
\newtheorem{proposition}[theorem]{Proposition~}
\newtheorem{corollary}[theorem]{Corollary~}
\newtheorem{definition}[theorem]{Definition~}
\newtheorem{notation}[theorem]{Notation~}
\newtheorem{example}[theorem]{Example~}
\newtheorem{problem}[theorem]{Problem~}
\newtheorem{remark}[theorem]{Remark~}

\newtheorem*{ac}{Acknowledgement}
\newtheorem{conjecture}{Conjecture~}

\newtheorem{main theorem}{Main Theorem}

\begin{document}
\title{Quantum smooth uncertainty principles for von Neumann bi-algebras} 
\author{Linzhe Huang}
\address{L. Huang, Yau Mathematical Sciences Center, Tsinghua University, Beijing, 100084, China}
\email{huanglinzhe@mail.tsinghua.edu.cn}

\author{Zhengwei Liu}
\address{Z. Liu, Yau Mathematical Sciences Center and Department of Mathematics, Tsinghua University, Beijing, 100084, China}
\address{Beijing Institute of Mathematical Sciences and Applications, Huairou District, Beijing, 101408, China}
\email{liuzhengwei@mail.tsinghua.edu.cn}

\author{Jinsong Wu}
\address{J. Wu, Institute for Advanced Study in Mathematics, Harbin Institute of Technology, Harbin, 150001, China}
\email{wjs@hit.edu.cn}
\date{}

\begin{abstract}
In this article, we prove various smooth uncertainty principles on von Neumann bi-algebras, which unify numbers of uncertainty principles on quantum symmetries, such as subfactors, and fusion bi-algebras etc, studied in quantum Fourier analysis. We also obtain Widgerson-Wigderson type uncertainty principles for von Neumann bi-algebras.
Moreover, we give a complete answer to a conjecture proposed by A. Wigderson and Y. Wigderson.
\end{abstract}

\maketitle

%
%
%



\section{Introduction}
Uncertainty principles have been investigated for more than hundred years in mathematics and physics inspired by the famous Heisenberg uncertainty principle
\cite{Hei,Ken, Weyl} with significant applications in information theory \cite{CRT,Don06}.

Recently quantum uncertainty principles on subfactors, an important type of quantum symmetires \cite{Jones0,EvaKaw}, have been established for support and for von Neumann entropy in \cite{JLW} and for R\'enyi entropy in \cite{LW}. These quantum uncertainty principles have been generalized on other types of quantum symmetries, such as Kac algebras \cite{LW2}, locally compact quantum groups \cite{JLW2} and fusion bialgebras \cite{LPW} etc, in the unified framework of quantum Fourier analysis \cite{JJLRW}. Such quantum inequalities were applied in the classification of subfactors \cite{Liu} and as analytic obstructions of unitary categorifications of fusion rings in \cite{LPW}.

In 2021, A. Wigderson and Y. Wigderson \cite{WW} introduced $k$-Hadamard matrices, as an analogue of discrete Fourier transforms, and they proved various uncertainty principles such as primary uncertainty principles, support uncertainty principles etc.
Their work unifies numbers of proofs of uncertainty principles in classical settings.

In this paper, we unify several quantum entropic uncertainty principles on quantum symmetries and we further generalize the results to various smooth entropies.
Inspired by the notion of $k$-Hadamard matrices, we introduce $k$-transforms between a pair of finite von Neumann algebras, and we call their combination a von Neumann $k$-bi-algebra.
We introduce various smooth entropies and prove the corresponding uncertainty principles for von Neumann $k$-bi-algebras.
On one hand, our results generalized numbers of uncertainty principles for quantum symmetries in \cite{JLW,LPW}. On the other hand, these results are slightly stronger
than uncertainty principles for $k$-Hadamard matrices in \cite{WW}. See Theorems \ref{thm:Quantum L1 smooth support  uncertainty principle}, \ref{thm:Quantum L2 smooth support  uncertainty principle}, \ref{thm:quantum_Hirschman_k_transform} and \ref{thm:quantum_smooth Hirschman_Beckner_k_transform}.

The primary uncertainty principle for $k$-Hadamard matrices plays a key role in \cite{WW} and we call this type of uncertainty principle the Wigderson-Wigderson uncertainty principle.
We prove the Wigderson-Wigderson uncertainty principle for von Neumann $k$-bi-algebras in Theorems \ref{thm:The  quantum  primary uncertainty principle} and for subfactors in Theorem \ref{thm:norm_quantum_Fourier_transform}.
In \cite{WW}, A. Wigderson and Y. Wigderson proposed a conjecture on the Wigderson-Wigderson uncertainty principle for the real line $\mathbb{R}$.
We give a complete answer to the conjecture, see Theorem \ref{thm:answer to conjecture} for details.

The paper is organized as follows.
In Section \ref{subsec:von Neumann bi-algebras}, we introduce $k$-transforms and von Neumann $k$-bi-algebras with examples from quantum Fourier analysis. We prove some basic uncertainty principles for von Neumann $k$-bi-algebras.
In Section \ref{sec:Quantum uncertainty principles}, we prove uncertainty principles on von Neumann bi-algebras for smooth support and von Neumann entropy perturbed by $p$-norms. We prove Wigderson-Wigderson uncertainty principles on von
Neumann bi-algebras, with a better constant in the case of subfactors.
In Section \ref{sec:conjecture}, we provide a bound for Wigderson-Wigderson uncertainty principle on the real line $\mathbb{R}$ and this answers a conjecture proposed by A. Wigderson and Y. Wigderson in \cite{WW}.

\begin{ac}
Zhengwei Liu was supported by NKPs (Grant no. 2020YFA0713000) and by Tsinghua University (Grant no. 100301004).
Jinsong Wu was supported by NSFC (Grant no. 11771413 and 12031004).
\end{ac}

\section{von Neumann bi-algebras and k-transforms}\label{subsec:von Neumann bi-algebras}
In this section, we recall some basic definitions and results about von Neumann algebras.
We introduce von Neumann bi-algebras with interesting examples and we prove some basic properties and uncertainty principles.

A von Neumann algebra $\mathcal{M}$ is said to be
\textbf{finite} if it has a faithful normal tracial positive linear functional $\tau_{\mathcal{M}}$, see e.g. \cite{KR}.
We will call this linear functional as trace in the rest of the paper.
We denote $\|x\|_p=\tau_{\mathcal{M}}(|x|^p)^{\frac{1}{p}}$, for $p>0$. When  $1\leq p<\infty$, $\|\cdot\|_p$ is called the $p$-norm.
Moreover, $\|x\|_{\infty}=\|x\|$, the operator norm of $x$. It is clear that $\|x\|_p=\|x^*\|_p=\||x|\|_p$ for $p>0$.

The following inequalities will be used frequently in the rest of the paper.
\begin{proposition}[H\"older's inequalities]\label{prop:Holder_inequality}
For any $x,y,z\in\mathcal{M}$, we have
\begin{enumerate}
\item $|\tau_{\mathcal{M}}(xy)|\leq\|x\|_p\|y\|_q$, where $1\leq p,q\leq\infty$, $\frac{1}{p}+\frac{1}{q}=1$;
\item $|\tau_{\mathcal{M}}(xyz)|\leq\|x\|_p\|y\|_q\|z\|_r$, where $1\leq p,q,r\leq\infty$, $\frac{1}{p}+\frac{1}{q}+\frac{1}{r}=1$;
\item $\|xy\|_r\leq\|x\|_p\|y\|_q$, where $0< p,q,r\leq\infty$, $\frac{1}{r}=\frac{1}{p}+\frac{1}{q}$.
\end{enumerate}
\begin{proof}
See e.g. Theorems 5.2.2 and 5.2.4 in \cite{Xu}.
\end{proof}

\end{proposition}

\begin{notation}
Suppose $\mathcal{A}$ and $\mathcal{B}$ are two finite von Neumann algebras with traces $d$ and $\tau$
respectively. Let $\mathcal{F}: \mathcal{A} \to \mathcal{B} $ be a linear map. For any $0< p,q \leq \infty$, define
$$
\|\mathcal{F}\|_{q\to p}:=\sup\{ \|\mathcal{F}(x)\|_p : x \in \mathcal{A},~ \|x\|_q=1 \}.
$$

\end{notation}

\begin{definition}\label{def:k_transform_von_bi-algebra}
Suppose $\mathcal{A}$ and $\mathcal{B}$ are two finite von Neumann algebras with traces $d$ and $\tau$
respectively. For $k>0$,
a \textbf{$k$-transform} $\mathcal{F}$ from $\mathcal{A}$ into $\mathcal{B}$ is a linear map such that $\|\mathcal{F}\|_{1\to\infty}\leq1$
and $\|\mathcal{F}^*\mathcal{F}(x)\|_{\infty}\geq k\|x\|_{\infty}$ for any $x\in\mathcal{A}$. We call the quintuple $(\mathcal{A},\mathcal{B},d,\tau,\mathcal{F})$ a \textbf{von Neumann $k$-bi-algebra}.
\end{definition}

\begin{example}\label{ex:Wigderson example Hadamard}
The definition of $k$-transform is inspired by the definition of $k$-Hadamard matrix of A. Wigderson and Y. Wigderson (Definition 2.2 in \cite{WW}). In particular, a $k$-Hadamard matrix $\mathcal{F}$ can be extended to a von Neumann $k$-bi-algebra  $(\mathcal{A},\mathcal{B},d,\tau,\mathcal{F})$, such that $\mathcal{A}$ and $\mathcal{B}$ are finite-dimensional abelian von Neumann algebras,  $d$ and $\tau$ are counting measures.
\end{example}

\begin{example}\label{ex:fushion bi-algebra}
Let the quintuple $(\mathcal{A},\mathcal{B},d,\tau,\mathcal{F})$ be a fusion bialgebra (See Definition 2.12 in \cite{LPW}),
where $\mathcal{A}$ and $\mathcal{B}$ are finite-dimensional von Neumann algebras with traces $d$ and $\tau$ respectively, and $\mathcal{A}$ is commutative, and $\mathcal{F}: \mathcal{A} \to \mathcal{B}$ is unitary with respect to $2$-norms.
By the quantum Hausdorff-Young inequality $\|\mathcal{F}\|_{1 \to \infty}=1$, (Theorem 4.5 in \cite{LPW}),
we have that $(\mathcal{A},\mathcal{B},d,\tau,\mathcal{F})$ is a von Neumann $1$-bi-algebra.

\end{example}

\begin{example}\label{ex:subfactor planar algebra}
Suppose $\mathscr{P}_{\bullet}$ is an irreducible subfactor planar algebra with finite Jones index (See Definition on page 4 in \cite{Jones0}) $\delta^2$,  $\delta>0$.
Let $Tr_{n,\pm}$ be the unnormalized Markov trace of $\mathscr{P}_{n,\pm}$, for $n\in \mathbb{N}$, and $\mathscr{F}: \mathscr{P}_{n,+} \to \mathscr{P}_{n,-}$ be the string Fourier transform, which is unitary.
Then by the quantum Hausdorff-Young inequality, (Theorem 4.8 and Theorem 7.3 in \cite{JLW}), we have that for any $n\in \mathbb{N}$, $2\leq p\leq\infty$ and $1/p+1/q=1$,
\begin{align*}
\|\mathscr{F}\|_{q \to p}=\bigg(\frac{1}{\delta}\bigg)^{1-\frac{2}{p}}.
\end{align*}
Therefore, $\|\delta \mathscr{F}\|_{1 \to \infty}=1$ and the quintuple $(\mathscr{P}_{n,+},\mathscr{P}_{n,-},Tr_{n,+},Tr_{n,-},\delta\mathscr{F})$ is a von Neumann $\delta^2$-bi-algebra.
\end{example}

\begin{remark}
The quantum Hausdorff-Young inequality, Theorem 7.3 in \cite{JLW}, also applies to reducible subfactor planar algebras, and in that case $\delta$ is replaced by certain constant $\delta_0$. Then $(\mathscr{P}_{n,+},\mathscr{P}_{n,-},Tr_{n,+},Tr_{n,-},\delta_0\mathscr{F})$ is a von Neumann $\delta_0^2$-bi-algebra.
\end{remark}

In \cite{WW}, Wigderson and Wigderson proved the primary uncertainty principles (See Theorem 2.3 in \cite{WW}) for any $k$-Hadamard matrix $A$,
\begin{align}\label{eq:wwup}
\|v\|_1\|Av\|_1\geq k\|v\|_{\infty}\|Av\|_{\infty},\quad v\in\mathbb{C}^n,
\end{align}
 which is the fundamental result of that paper. We call the inequality as Wigderson-Wigderson uncertainty principle. In this paper, we prove the following quantum version of Wigderson-Wigderson uncertainty principle for von Neumann $k$-bi-algebras.
 When a von Neumann $k$-bi-algebra is obtained from Example \ref{ex:Wigderson example Hadamard}, then our theorem implies Theorem 2.3 in \cite{WW}.
\begin{theorem}[The  quantum Wigderson-Wigderson uncertainty principle]\label{thm:The  quantum  primary uncertainty principle}
Let $(\mathcal{A},\mathcal{B},d,\tau,\mathcal{F})$
be a von Neumann $k$-bi-algebra.
For any $x\in\mathcal{A}$, we have
\[\|x\|_1\|\mathcal{F}(x)\|_1\geq k\|x\|_{\infty}\|\mathcal{F}(x)\|_{\infty}.\]
\begin{proof}
When $1\leq p, q \leq\infty$ and $1/p+1/q=1$, we have that $\|\mathcal{F}^*\|_{p \to q}=\|\mathcal{F}\|_{p \to q}$, because
\begin{align*}
\|\mathcal{F}\|_{p\to q}=&\sup\{ \|\mathcal{F}(x) \|_q : x\in \mathcal{A}, \|x\|_p=1 \} \\
=&\sup\{ |\tau(\mathcal{F}(x) y^* ) | : x\in \mathcal{A}, y \in \mathcal{B}, \|x\|_p=1, \|y\|_p=1\} \\
=&\sup\{  |d( x (\mathcal{F}^*(y) )^* ) | :  x\in \mathcal{A}, y \in \mathcal{B}, \|x\|_p=1, \|y\|_p=1\} \\
=&\sup\{  \| \mathcal{F}^*(y)  \|_q : y \in \mathcal{B}, \|y\|_p=1\} \\
=&\|\mathcal{F}^*\|_{p\to q} \;.
\end{align*}
This implies that  $\|\mathcal{F}^*\|_{1\to\infty}=\|\mathcal{F}\|_{1\to\infty}\leq1$.
Then for any $x\in\mathcal{A}$, we have
\begin{align*}
\|\mathcal{F}(x)\|_{\infty}\leq\|x\|_1,\quad k\|x\|_{\infty}\leq\|\mathcal{F}^*\mathcal{F}(x)\|_{\infty}\leq\|\mathcal{F}(x)\|_1.
\end{align*}
 Multiplying the above two inequalities, we obtain
\begin{align*}
\|x\|_1\|\mathcal{F}(x)\|_1\geq k\|x\|_{\infty}\|\mathcal{F}(x)\|_{\infty}.
\end{align*}
This completes the proof of the theorem.
\end{proof}
\end{theorem}
Using the primary uncertainty principle, A. Wigderson and Y. Wigderson further prove the Donoho-Stark uncertainty principle for arbitrary $k$-Hadamard matrices (See Theorem 3.2 in \cite{WW}). In this paper, we prove the
Donoho-Stark uncertainty principle for von Neumann $k$-bi-algebras using the quantum Wigderson-Wigderson uncertainty principle.
Firstly, let's recall the notion of the support in a finite von Neumann algebra.

\begin{definition}
Let $\mathcal{M}$ be a finite von Neumann algebra with a trace $\tau_{\mathcal{M}}$.
For any $x\in\mathcal{M}$, let $\mathcal{R}(x)$ be the range projection of $x$. The \textbf{support $\mathcal{S}(x)$} of $x$
is defined as $\tau_{\mathcal{M}}(\mathcal{R}(x))$.
\end{definition}

The support has been used in the  quantum Donoho-Stark uncertainty principles on quantum symmetries such as subfactors and fusion rings, see Theorem 5.2 in \cite{JLW} and Theorem 4.8 in \cite{LPW} respectively. We generalize the Donoho-Stark uncertainty principles from these quantum symmetries to von Neumann $k$-bi-algebras.

\begin{theorem}[Quantum Donoho-Stark  uncertainty principle]\label{thm:quantum_support-size _uncertainty_principle}
Let $(\mathcal{A},\mathcal{B},d,\tau,\mathcal{F})$
be a von Neuman $k$-bi-algebra.
Then for any non-zero operator  $x\in\mathcal{A}$, we have
\[\mathcal{S}(x)\mathcal{S}(\mathcal{F}(x))\geq k.\]
\begin{proof}
We already have, from Theorem \ref{thm:The  quantum  primary uncertainty principle}, that for any nonzero $x\in\mathcal{A}$,
\begin{align*}
\|x\|_1\|\mathcal{F}(x)\|_1\geq k\|x\|_{\infty}\|\mathcal{F}(x)\|_{\infty}.
\end{align*}
Thus, all we need is to bound the 1-norm by the support of $x$, which can be implemented through H\"older's inequality, for any  $x\in \mathcal{A}$,
\begin{align*}
\|x\|_1=\|\mathcal{R}(x)x\|_1\leq\|x\|_{\infty}\|\mathcal{R}(x)\|_1=\|x\|_{\infty}\mathcal{S}(x).
\end{align*}
Applying this bound to both $x$ and $\mathcal{F}(x)$, we obtain the result.

\end{proof}
\end{theorem}
\begin{remark}
Our theorem is a generalization of the Donoho-Stark uncertainty principle in \cite{DS} and some variations,
\begin{enumerate}
\item In Example \ref{ex:Wigderson example Hadamard}, Theorem \ref{thm:quantum_support-size _uncertainty_principle} implies Theorem 3.2 in \cite{WW};
\item In Example \ref{ex:fushion bi-algebra},  Theorem \ref{thm:quantum_support-size _uncertainty_principle} implies Theorem 4.8 in \cite{LPW};
\item In Example \ref{ex:subfactor planar algebra},  Theorem \ref{thm:quantum_support-size _uncertainty_principle} implies Theorem 5.2 in \cite{JLW}.
\end{enumerate}
\end{remark}

\section{Quantum smooth uncertainty principles}
\label{sec:Quantum uncertainty principles}
In this section, we prove a series of smooth uncertainty principles for von Neumann bi-algebras. We firstly prove the quantum smooth support uncertainty principles in \S \ref{subsec:Quantum smooth support uncertainty principles}. Then we proceed to prove quantum Wigderson-Wigderson uncertainty principles for general $p$-norms, $1\leq p \leq \infty$, and give an example concerning the quantum Fourier transform on subfactor planar algebras in \S \ref{subsec:Quantum Wigderson Wigderson uncertainty principle}.
Finally, we also prove quantum smooth Hirschman-Becker uncertainty principles in \S \ref{subsec:Quantum_ Hirschman_Beckner_uncertainty_principle}.

\subsection{Quantum smooth support uncertainty principles}
\label{subsec:Quantum smooth support uncertainty principles}We firstly introduce a new smooth support which is slightly different from the classical smooth support.
\begin{definition}\label{def:quantum smooth support}
Let $\mathcal{M}$ be a finite von Neumann algebera with a trace $\tau_{\mathcal{M}}$.
Let $\epsilon\in [0,1]$ and $p\in [1,\infty]$. For any element $x\in\mathcal{M}$,
we define the $(p,\epsilon)$  smooth support to be
\begin{align*}
\mathcal{S}_{\epsilon}^p(x)=\inf\{\tau_{\mathcal{M}}(H\mathcal{R}(x)):
H\in \mathcal{M},\ 0\leq H\leq I,\ \|(I-H)x\|_p\leq \epsilon \|x\|_p\},
\end{align*}
where $\mathcal{R}(x)$ is the range projection of $x$.
\end{definition}

\begin{remark}
Since the set
$$\mathscr{S}(\epsilon,p,x):=\{H\in \mathcal{M}:\ \ 0\leq H\leq I,\ \|(I-H)x\|_p\leq\epsilon\|x\|_p\}$$
is compact in the weak operator topology and the trace is normal,
there exits an $H_0\in\mathscr{S}$ such that $\mathcal{S}_{\epsilon}^p(x)=\tau_{\mathcal{M}}(H_0\mathcal{R}(x))$.
\end{remark}

\begin{remark}
Take $\epsilon=0$, then $(I-H)x=0$ and this implies
$H\mathcal{R}(x)=\mathcal{R}(x)$. In this case, $\mathcal{S}_0^p(x)=\mathcal{S}(x)$.
\end{remark}

Besides Definition \ref{def:quantum smooth support}, there are three kinds of notions of the smooth support.

\begin{definition}\label{def: three support}
Let $\mathcal{M}$ be a finite von Neumann algebra with a trace $\tau_{\mathcal{M}}$.
Let $\epsilon\in [0,1]$ and $p\in [1,\infty]$. For any element $x\in\mathcal{M}$,
define
\begin{align*}
f_1(\epsilon,p,x)&:=\inf\{\tau_{\mathcal{M}}(\mathcal{R}(y)):
y\in \mathcal{M},\ \|x-y\|_p\leq \epsilon \|x\|_p\}, \\
f_2(\epsilon,p,x)&:=\inf\{\tau_{\mathcal{M}}(\mathcal{R}(Hx)):
H\in \mathcal{M},\  0\leq H\leq I,\ \|(I-H)x\|_p\leq \epsilon \|x\|_p\},\\
f_3(\epsilon,p,x)&:=\inf\{\tau_{\mathcal{M}}(\mathcal{R}(Qx)):
Q\in \mathcal{M},\ Q=Q^*=Q^2,\ \|(I-Q)x\|_p\leq \epsilon \|x\|_p\}.
\end{align*}

\end{definition}

\begin{proposition}\label{prop:three support equivalent and our support}
For any $x\in\mathcal{M}$, we have $$ \mathcal{S}_{\epsilon}^p(x)\leq f_1(\epsilon,p,x)=f_2(\epsilon,p,x)=f_3(\epsilon,p,x).$$
\begin{proof}
It is clear that  $f_1(\epsilon,p,x)\leq f_2(\epsilon,p,x)\leq f_3(\epsilon,p,x)$.

For any $y\in\mathcal{M}$, we claim that
$$\tau_{\mathcal{M}}(\mathcal{R}(\mathcal{R}(y)x))\leq\tau_{\mathcal{M}}(\mathcal{R}(y)),\quad \|(I-\mathcal{R}(y))x\|_p\leq\|x-y\|_p.$$
If the claim holds, then $f_3(\epsilon,p,x)\leq f_1(\epsilon,p,x)$.
Since $\mathcal{R}(\mathcal{R}(y)x)\leq \mathcal{R}(y)$, the first inequality holds.

Next, we prove the second inequality in the claim.
It is enough to prove that $|x-\mathcal{R}(y)x|\leq |x-y|$. Since $\sqrt{\cdot}$ is an operator-monotone function, it reduces to prove $(x-\mathcal{R}(y)x)^*(x-\mathcal{R}(y)x)\leq (x-y)^*(x-y)$.
For any normal state $\rho$ on $\mathcal{M}$,
by the Cauchy-Schwartz inequality, we have
\begin{align*}
2|\rho(y^*x)|
=&2|\left\langle x,y\right\rangle_{\rho}|\\
=& 2|\left\langle\mathcal{R}(y)x,y\right\rangle_{\rho}| \\
\leq & 2\left\langle \mathcal{R}(y)x,\mathcal{R}(y)x \right\rangle_{\rho}^{\frac{1}{2}}\left\langle y,y\right\rangle_{\rho}^{\frac{1}{2}} \\
\leq & \rho(x^*\mathcal{R}(y)x)+\rho(y^*y).
\end{align*}
Therefore,
\begin{align*}
\rho(x^*y)+\rho(y^*x)\leq
\rho(x^*\mathcal{R}(y)x)+\rho(y^*y).
\end{align*}
Rearranging the above inequality, we obtain
\begin{align*}
\rho((x-\mathcal{R}(y)x)^*(x-\mathcal{R}(y)x))\leq
\rho((x-y)^*(x-y)).
\end{align*}
Thus,
$$(x-\mathcal{R}(y)x)^*(x-\mathcal{R}(y)x)\leq (x-y)^*(x-y).$$
The claim holds and we have $f_1(\epsilon,p,x)=f_2(\epsilon,p,x)=f_3(\epsilon,p,x)$.

For any $H\in\mathcal{M}$, $0\leq H\leq I$,
we have
\begin{align*}
\tau_{\mathcal{M}}(\mathcal{R}(x)H)\leq\tau_{\mathcal{M}}(|\mathcal{R}(x)H|)\leq
\tau_{\mathcal{M}}(\|\mathcal{R}(x)H\|\mathcal{R}(Hx))\leq\tau_{\mathcal{M}}(\mathcal{R}(Hx)).
\end{align*}
The first inequality is true by H\"older's inequality.
The second one uses the fact that $|y^*|\leq\|y\|\mathcal{R}(y)$, $y\in\mathcal{M}$.
The last inequality is due to $\|\mathcal{R}(x)H\|\leq1$.
So we have
$$\mathcal{S}_{\epsilon}^p(x)\leq f_2(\epsilon,p,x).$$
In summary, the statement holds.
\end{proof}
\end{proposition}
In \cite{WW}, A. Wigderson and Y. Wigderson introduced the following smooth support for the finite-dimensional and abelian case.
\begin{definition}(See Definition 3.15 in \cite{WW})
\label{def: smooth support Wigderson}
Let $\mathcal{M}=\mathbb{C}^n$, $n\in\mathbb{N}^*$, and $\tau_{\mathcal{M}}$ be the counting measure.
Let $\epsilon\in [0,1]$ and $p\in [1,\infty]$. For an operator $x\in\mathcal{M}$, the $(p,\epsilon)$ support-size of $x$ is defined to be
\begin{align*}
|supp_{\epsilon}^p(x)|=\min\{\tau_{\mathcal{M}}(Q):\
Q\in\mathcal{M},\ Q=Q^*=Q^2,\
 \|(I-Q)x\|_p\leq\epsilon\|x\|_p\}.
\end{align*}
\end{definition}

\begin{remark}\label{rem: our support bound by WW support}
When $\mathcal{M}$ is finite-dimensional and abelian and $\tau_{\mathcal{M}}$ is the counting measure, then $f_3(\epsilon,p,x)$ is equal to $|supp_{\epsilon}^p(x)|$. In this case,  $\mathcal{S}_{\epsilon}^p(x)\leq|supp_{\epsilon}^p(x)|$.
\end{remark}

\begin{lemma}For any $x\in\mathcal{M}$, we have
$\mathcal{S}_{\epsilon}^p(x)$ is continuous with respect to $\epsilon$.
\begin{proof}
When $0<c<1$, take an $H\in\mathcal{M}$ such that
\begin{align*}
\mathcal{S}_{\epsilon}^p(x)=\tau_{\mathcal{M}}(H\mathcal{R}(x)),\quad \|(I-H)x\|_p\leq\epsilon\|x\|_p,\quad 0\leq H\leq I.
\end{align*}
Let $H'=I-c(I-H)$, then $0\leq H'\leq I$.
Moreover, we have
\begin{align*}
\tau_{\mathcal{M}}(H'\mathcal{R}(x))=(1-c)\tau_{\mathcal{M}}(\mathcal{R}(x))+c \mathcal{S}_{\epsilon}^p(x),\quad \|(I-H')x\|_p\leq c\epsilon\|x\|_p.
\end{align*}
Therefore,
\begin{align}\label{eq:epsilon c}
\mathcal{S}_{\epsilon}^p(x)\leq\mathcal{S}_{c\epsilon}^p(x)\leq (1-c)\tau_{\mathcal{M}}(\mathcal{R}(x))+c \mathcal{S}_{\epsilon}^p(x).
\end{align}
So \[\lim_{c\to1^-}\mathcal{S}_{c\epsilon}^p(x)=\mathcal{S}_{\epsilon}^p(x).\]

When $c>1$, replacing $c$ by $c^{-1}$ and $\epsilon$
by $c\epsilon$ in Inequality (\ref{eq:epsilon c}), we have
$$\mathcal{S}_{c\epsilon}^p(x)\leq\mathcal{S}_{\epsilon}^p(x)\leq \left(1-\frac{1}{c}\right)\tau_{\mathcal{M}}(\mathcal{R}(x))+\frac{1}{c} \mathcal{S}_{c\epsilon}^p(x).$$
So
 \[\lim_{c\to1^+}\mathcal{S}_{c\epsilon}^p(x)=\mathcal{S}_{\epsilon}^p(x).\]
From the above discussions, $\mathcal{S}_{\epsilon}^p(x)$ is continuous with respect to $\epsilon$.

\end{proof}
\end{lemma}

We have the following
quantum $L^1$ smooth support uncertainty principle.
\begin{theorem}[The quantum $L^1$ smooth support  uncertainty principle]\label{thm:Quantum L1 smooth support  uncertainty principle}
Let the quintuple
 $(\mathcal{A},\mathcal{B},d,\tau, \mathcal{F})$
be a von Neumann $k$-bi-algebra and $x\in\mathcal{A}$ be a non-zero operator. For any $\epsilon,\eta\in [0,1]$, we have
\begin{align*}
\mathcal{S}_{\epsilon}^1(x)\mathcal{S}_{\eta}^1(\mathcal{F}(x))\geq k(1-\epsilon)(1-\eta).
\end{align*}
\begin{proof}
Take a  positive operator $H$ in $\mathcal{A}$
such that
$$\mathcal{S}_{\epsilon}^1(x)=d(\mathcal{R}(x)H),\quad \|(I-H)x\|_1\leq\epsilon\|x\|_1,\quad
0\leq H\leq I.$$
By H\"older's inequality, we have
\begin{align*}
d(|x^*|(I-H))\leq d(|x^*(I-H)|)=\|(I-H)x\|_1\leq\epsilon\|x\|_1.
\end{align*}
 Thus
\begin{align*}
\mathcal{S}_{\epsilon}^1(x)&=d(\mathcal{R}(x)H)\\
&=\frac{1}{\|x\|_{\infty}}d(\|x\|_{\infty}\mathcal{R}(x)H)\\
&\geq \frac{1}{\|x\|_{\infty}}d(|x^*|H)=\frac{1}{\|x\|_{\infty}}d(|x^*|)-\frac{1}{\|x\|_{\infty}}d(|x^*|(I-H))\\
&\geq\frac{\|x\|_1}{\|x\|_{\infty}}-\frac{1}{\|x\|_{\infty}}d(|x^*(I-H)|)\\
&\geq\frac{\|x\|_1}{\|x\|_{\infty}}(1-\epsilon).
\end{align*}
Repeating the above process for $\mathcal{F}(x)$, we obtain
\begin{align*}
\mathcal{S}_{\eta}^1(\mathcal{F}(x))\geq
\frac{\|\mathcal{F}(x)\|_1}{\|\mathcal{F}(x)\|_{\infty}}(1-\eta).
\end{align*}
Multiplying these two inequalities, we have
\begin{align*}
\mathcal{S}_{\epsilon}^1(x)\mathcal{S}_{\eta}^1(\mathcal{F}(x))\geq \frac{\|x\|_1}{\|x\|_{\infty}}\cdot
\frac{\|\mathcal{F}(x)\|_1}{\|\mathcal{F}(x)\|_{\infty}}(1-\epsilon)(1-\eta)
\geq k(1-\epsilon)(1-\eta).
\end{align*}
The second inequality uses Theorem \ref{thm:The  quantum  primary uncertainty principle}, the quantum Wigderson-Wigderson uncertainty principle.
\end{proof}
\end{theorem}

\begin{remark}
We can obtain Theorem \ref{thm:quantum_support-size _uncertainty_principle} from Theorem \ref{thm:Quantum L1 smooth support  uncertainty principle} by assuming $\epsilon=\eta=0$.
\end{remark}

Applying Theorem \ref{thm:Quantum L1 smooth support  uncertainty principle} to the quantum Fourier transform on subfactor planar algebras, we obtain the following corollary.

\begin{corollary}\label{cor:L1_quantum_Fourier}
Suppose $\mathscr{P}_{\bullet}$ is an irreducible subfactor planar algebra with finite Jones index $\delta^2$. Let $\mathscr{F}$ be the Fourier transform from $\mathscr{P}_{n,\pm}$ onto $\mathscr{P}_{n,\mp}$.
Then for any non-zero $n$-box $x\in\mathscr{P}_{n,\pm}$, we have
\begin{align*}
\mathcal{S}_{\epsilon}^1(x)\mathcal{S}_{\eta}^1(\mathscr{F}(x))\geq \delta^2(1-\epsilon)(1-\eta).
\end{align*}
\end{corollary}

When $p=2$ in Definition \ref{def:quantum smooth support}, we are able to choose a positive contraction $H$ in the abelian *-subalgebra generated by $|x^*|$ such that the $(2,\epsilon)$ support-size is exactly the trace of $H$. More precisely, we have

\begin{proposition}\label{prop:p=2}Suppose $\mathcal{M}$ is a finite von Neumann algebra with a  trace $\tau_{\mathcal{M}}$.
Let $x\in\mathcal{M}$, and let $\mathcal{N}$ be the abelian von Neumann subalgebra
generated by $|x^*|$ in $\mathcal{M}$.
For any $\epsilon\in[0,1]$, we have
\begin{align*}
\mathcal{S}_{\epsilon}^2(x)=\min\{\tau_{\mathcal{M}}(H):
H\in \mathcal{N},\ 0\leq H\leq \mathcal{R}(x),\ \|(I-H)x\|_2\leq \epsilon \|x\|_2\}.
\end{align*}

\begin{proof}

Let $\Phi$ be the trace-preserving conditional expectation from $\mathcal{M}$ onto $\mathcal{N}$.
For any $H \in \mathcal{M}$, $0\leq H \leq I$. Take $H'=\Phi(H)\mathcal{R}(x)$, then
$$\tau_{\mathcal{M}}(H')= \tau_{\mathcal{M}}(\Phi(H)\mathcal{R}(x)) =\tau_{\mathcal{M}}(\mathcal{R}(x)H),$$
and
$H' \in \mathbb{N}$ and $0\leq H' \leq \mathcal{R}(x)$.

Note that any pure state $\rho$ on $\mathcal{N}$ is multiplicative, so $\rho(|\Phi(y)|^2)=|\rho\circ\Phi(y)|^2$, for any $y \in \mathcal{M}$.
Moreover. $\rho\circ\Phi$ is a state on $\mathcal{M}$, by the Cauchy-Schwartz inequality,
$|\rho\circ\Phi(y)|^2\leq \rho \circ \Phi (|y|^2)$.
So $\rho(|\Phi(y)|^2) \leq \rho(\Phi(|y|^2))$, and therefore
$|\Phi(y)|^2\leq \Phi (|y|^2)$.

Take $y=I-H$, then
\begin{align*}
\|(I-H')x\|_2^2&=\|\Phi(I-H)\mathcal{R}(x)x\|_2^2\\
&=\tau_{\mathcal{M}}(|\Phi(I-H)|^2|x^*|^2)\\
&\leq \tau_{\mathcal{M}}(\Phi(|I-H|^2)|x^*|^2)\\
&=\tau_{\mathcal{M}}(|I-H|^2|x^*|^2)\\
&=\|(I-H)x\|_2^2.
\end{align*}
Therefore, the statement holds.

\end{proof}

\color{black}

\end{proposition}
We have the following
quantum $L^2$  smooth support  uncertainty principle.

\begin{theorem}[The quantum $L^2$ smooth support  uncertainty principle]\label{thm:Quantum L2 smooth support  uncertainty principle}
Let $(\mathcal{A},\mathcal{B},d,\tau,\mathcal{F})$
be a von Neumann $k$-bi-algebra.
Suppose $\mathcal{A}$ and $\mathcal{B}$ are finite dimensional and $\mathcal{F}^*\mathcal{F}=kI$.
For any non-zero operator $x\in\mathcal{A}$, we have
\begin{align*}
\mathcal{S}_{\epsilon}^2(x)\mathcal{S}_{\eta}^2(\mathcal{F}(x))\geq k(1-\epsilon-\eta)^2,
\quad \forall\epsilon,\eta\in [0,1],\ \epsilon+\eta\leq1.
\end{align*}
\end{theorem}
\begin{proof}
Take $W=\mathcal{F}/\sqrt{k}$, then $W^*W=I$.
Since the definition of $\mathcal{S}^2_\eta$ is invariant under rescaling, we have that $\mathcal{S}_{\eta}^2(W(x))=\mathcal{S}_{\eta}^2(\mathcal{F}(x))$.

Let $x=|x^*|U$ and $y=W(x)=|y^*|V$ be the polar decompositions, where $U$ and $V$ are the polar parts in
$\mathcal{A}$ and $\mathcal{B}$ respectively. Let $\mathcal{A}_0$ be the abelian von Neumann subalgebra of $\mathcal{A}$ generated by $|x^*|$ and $\mathcal{B}_0$ be the abelian von Neumann subalgebra of $\mathcal{B}$ generated by $|y^*|$.
Let $\Phi$ be the trace-preserving conditional expectation from $\mathcal{B}$ onto $\mathcal{B}_0$ and $M=\Phi\mathcal{R}_{ V^*}W\mathcal{R}_{U}$. Then $M$ is a linear operator from
 $\mathcal{A}_0$ into  $\mathcal{B}_0$ such that $M|x^*|=|y^*|$. Let
$\{e_i\}_{i=1}^n$ and $\{f_j\}_{j=1}^m$
be mutually orthogonal minimal projections in $\mathcal{A}_{0}$ and $\mathcal{B}_{0}$ such that $\sum_{i=1}^n e_i=I_{\mathcal{A}}$ and $\sum_{j=1}^m f_j=I_{\mathcal{B}}$.
The linear operator $M$ is a $m\times n$ matrix $(a_{ij})_{1\leq i\leq m,1\leq j\leq n}$ with $|a_{ij}|\leq d(e_j)/\sqrt{k}$ by $\|W\|_{1\to\infty}\leq1/\sqrt{k}$.

By Proposition \ref{prop:p=2}, we can
find two  positive operators $H$ in $\mathcal{A}_0$ and $K$ in $\mathcal{B}_0$ such that
\begin{align*}
0\leq H\leq \mathcal{R}(x)&,\quad 0\leq K\leq \mathcal{R}(y),\\
\|(I-H)x\|_2\leq\epsilon\|x\|_2&,\quad\|(I-K)y\|_2\leq\eta\|y\|_2,\\
d(H)=\mathcal{S}_{\epsilon}^2(x)&,\quad
\tau(K)=\mathcal{S}_{\eta}^2(y).
\end{align*}
  By direct computations, we have
\begin{align*}
H=\sum_{i=1}^n\frac{d(e_i H)}{d(e_i)}e_i,&\quad
 K=\sum_{j=1}^m\frac{\tau(f_j K)}{\tau(f_j)}f_j.
\end{align*}
Let $\widetilde{M}=KMH$, then $\widetilde{M}$ is a linear operator from $\mathcal{A}_0$ into $\mathcal{B}_0$.
For any $v\in\mathcal{A}_0$, we have
\begin{align*}
\|\widetilde{M}v\|_2^2&=\sum_{i=1}^m \tau(f_i)\bigg|\frac{\tau(f_i K)}{\tau(f_i)}\sum_{j=1}^na_{ij}\frac{d(e_j H)}{d(e_j)}v_j\bigg|^2\\
&\leq \sum_{i=1}^m
\frac{\tau(f_i K)^2}{\tau(f_i)}\sum_{j=1}^n|a_{ij}|^2\frac{d(e_j H)^2}{d(e_j)^3}\sum_{j=1}^n d(e_j)|v_j|^2\\
&=\sum_{i=1}^m
\frac{\tau(f_i K)^2}{\tau(f_i)}\sum_{j=1}^n|a_{ij}|^2\frac{d(e_j H)^2}{d(e_j)^3}\|v\|_2^2\\
&\leq\frac{1}{k}\sum_{i=1}^m
\frac{\tau(f_i K)^2}{\tau(f_i)}\sum_{j=1}^n\frac{d(e_j H)^2}{d(e_j)}\|v\|_2^2\\
&\leq \frac{d(H)\tau(K)}{k}\|v\|_2^2.
\end{align*}
The first inequality is true by the Cauchy-Schwartz inequality and the second one uses the fact that $|a_{ij}|\leq d(e_j)/\sqrt{k}$.
This implies
\begin{align}\label{eq:upper bound}
\|\widetilde{M}\|_{2\to2}\leq\sqrt{d(H)\tau(K)}/\sqrt{k}=\sqrt{\mathcal{S}_{\epsilon}^2(x)\mathcal{S}_{\eta}^2(y)}/\sqrt{k}.
\end{align}

For the lower bound of $\widetilde{M}$, we firstly observe that
\begin{align*}
\|M(I-H)|x^*|\|_2&=\|\Phi\mathcal{R}_{V^*}W\mathcal{R}_U(I-H)|x^*|\|_2\\
&\leq
\|\mathcal{R}_{V^*}W\mathcal{R}_U(I-H)|x^*|\|_2 \\
&=\|(I-H)|x^*|\|_2.
\end{align*}
Since $K$ is a contraction, so
\begin{align*}
\|KM(I-H)|x^*|\|_2\leq\|M(I-H)|x^*|\|_2\leq\|(I-H)|x^*|\|_2.
\end{align*}
Therefore, we have
\begin{align*}
\|M|x^*|-\widetilde{M}|x^*| \|_2&=\|M|x^*|-KMH|x^*|\|_2 \\
&=\|(I-K)M|x^*|+ KM(I-H)|x^*|\|_2\\
&\leq\|(I-K)|y^*|\|_2+\|(I-H)|x^*|\|_2 \\
& \leq (\epsilon+\eta)\||x^*|\|_2.
\end{align*}
This implies
\begin{align}\label{eq:lower bound}
\begin{split}
\|\widetilde{M}|x^*|\|_2&\geq\|M|x^*|\|_2-(\epsilon+\eta)\||x^*|\|_2\\
&=\||y^*|\|_2-(\epsilon+\eta)\||x^*|\|_2\\
& =(1-\epsilon-\eta)\||x^*| \|_2.
\end{split}
\end{align}
Finally, combining equations (\ref{eq:upper bound}) and (\ref{eq:lower bound}) we see that
\begin{align*}
\mathcal{S}_{\epsilon}^2(x)\mathcal{S}_{\eta}^2(\mathcal{F}(x))\geq k(1-\epsilon-\eta)^2.
\end{align*}
This completes the proof of the theorem.
\end{proof}

\begin{remark}\label{Rem: Supp UP}
When $\mathcal{F}$ is a $k$-Hadamard matrix,
A. Wigderson and Y. Wigderson proved the following results (See Theorems 3.17 and 3.20 in \cite{WW} ):
\begin{enumerate}
\item For any $x\in\mathcal{M}$,
\[  |supp_{\epsilon}^1(x)||supp_{\eta}^1(\mathcal{F}(x))|\geq (1-\epsilon)(1-\eta),\quad \forall\epsilon,\eta\in[0,1]; \]
\item If $\mathcal{F}^*\mathcal{F}=kI$, then
for any $x\in\mathcal{M}$,
\[ |supp_{\epsilon}^2(x)||supp_{\eta}^2(\mathcal{F}(x))|\geq (1-\epsilon-\eta)^2,\quad\forall \epsilon,\eta\in[0,1],\quad\epsilon+\eta\leq1.\]
\end{enumerate}
By Remark \ref{rem: our support bound by WW support}, we have
\begin{align*}
 |supp_{\epsilon}^1(x)||supp_{\eta}^1(\mathcal{F}(x))|&\geq\mathcal{S}_{\epsilon}^1(x)\mathcal{S}_{\eta}^1(\mathcal{F}(x))\geq (1-\epsilon)(1-\eta),\\
 |supp_{\epsilon}^2(x)||supp_{\eta}^2(\mathcal{F}(x))|&\geq\mathcal{S}_{\epsilon}^2(x)\mathcal{S}_{\eta}^2(\mathcal{F}(x))\geq (1-\epsilon-\eta)^2.
\end{align*}
So Theorems \ref{thm:Quantum L1 smooth support  uncertainty principle} and \ref{thm:Quantum L2 smooth support  uncertainty principle} imply Theorems 3.17 and 3.20 in \cite{WW}.
\end{remark}
When $\mathcal{F}$ is a $k$-Hadamard matrix, Theorems \ref{thm:Quantum L1 smooth support  uncertainty principle} and \ref{thm:Quantum L2 smooth support  uncertainty principle} are strictly stronger than
 Theorems 3.17 and 3.20 in \cite{WW}. We construct the following example.
\begin{example}
Let $\mathcal{A}=\mathcal{B}=\mathbb{C}\oplus\mathbb{C}$ and $d(f)=\tau(f)=f(0)+f(1)$, $f\in \mathbb{C}^2$.
Take $x=(1,1)\in \mathbb{C}^2$ and $\epsilon=\eta=1/3$. Then $|supp_{\epsilon}^1(x)|=|supp_{\epsilon}^2(x)|=2$ while $\mathcal{S}_{\epsilon}^1(x)=\mathcal{S}_{\epsilon}^2(x)=4/3$. Let $\mathcal{F}=I$ be the 1-transform, we have
\begin{align*}
4&=|supp_{\epsilon}^1(x)||supp_{\eta}^1(\mathcal{F}(x))|>
\mathcal{S}_{\epsilon}^1(x)\mathcal{S}_{\eta}^1(\mathcal{F}(x))=\frac{16}{9},\\
4&=|supp_{\epsilon}^2(x)||supp_{\eta}^2(\mathcal{F}(x))|>
\mathcal{S}_{\epsilon}^2(x)\mathcal{S}_{\eta}^2(\mathcal{F}(x))=\frac{16}{9}.
\end{align*}
\end{example}
Applying Theorem \ref{thm:Quantum L2 smooth support  uncertainty principle} to the quantum Fourier transform on subfactor planar algebras, we obtain the following corollary.

\begin{corollary}\label{cor:L2_quantum_Fourier}
Suppose $\mathscr{P}_{\bullet}$ is an irreducible subfactor planar algebra with finite Jones index $\delta^2$. Let $\mathscr{F}$ be the Fourier transform from $\mathscr{P}_{n,\pm}$ onto $\mathscr{P}_{n,\mp}$.
Then for any non-zero $n$-box $x\in\mathscr{P}_{n,\pm}$, we have
\begin{align*}
\mathcal{S}_{\epsilon}^2(x)\mathcal{S}_{\eta}^2(\mathcal{F}(x))\geq k(1-\epsilon-\eta)^2,
\quad \forall\epsilon,\eta\in [0,1],\ \epsilon+\eta\leq1.
\end{align*}
\end{corollary}

\subsection{Quantum Wigderson-Wigderson uncertainty principle}
\label{subsec:Quantum Wigderson Wigderson uncertainty principle}
In this section, we prove the quantum Wigderson-Wigderson uncertainty principle for von Neumann $k$-bi-algebras for $1/p+1/q=1$, and for quantum Fourier transform on subfactor planar algebras for any $0<p,q\leq\infty$.

We prove the quantum Hausdorff-Young inequality for $k$-transforms using the standard interpolation method.
\begin{theorem}\label{thm:Hausdorff_Young_k_transform}
Let $(\mathcal{A},\mathcal{B},d,\tau,\mathcal{F})$
be a von Neumann $k$-bi-algebra such that $\mathcal{F}^*\mathcal{F}=kI$. For any $x\in\mathcal{A}$, we have
\begin{align*}
\|\mathcal{F}(x)\|_p\leq k^{\frac{1}{p}}\|x\|_q,
\end{align*}
where $2\leq p\leq\infty$ and $1/p+1/q=1$.
\begin{proof}
Note that
\begin{align*}
\|\mathcal{F}(x)\|_{\infty}\leq\|x\|_1,\quad\|\mathcal{F}(x)\|_2=\sqrt{k}\|x\|_2.
\end{align*}
Applying the Riesz-Thorin interpolation theorem (\cite{Ree}, Theorem IX.17),
we have that $\|\mathcal{F}(x)\|_p\leq k^{\frac{1}{p}}\|x\|_q$.
\end{proof}
\end{theorem}

Then we have the following quantum Wigderson-Wigderson uncertainty principles for $k$-transforms.
\begin{theorem}\label{thm:quantum_Wigderson_uncertainty_principle_k_transform}
Let $(\mathcal{A},\mathcal{B},d,\tau,\mathcal{F})$
be a von Neumann $k$-bi-algebra such that $\mathcal{F}^*\mathcal{F}=kI$. For any $x\in\mathcal{A}$, we have
\begin{align*}
\|x\|_q\|\mathcal{F}(x)\|_q\geq k^{1-\frac{2}{p}}\|x\|_p\|\mathcal{F}(x)\|_p,
\end{align*}
where $2\leq p\leq\infty$ and $1/p+1/q=1$.
\end{theorem}
\begin{proof}
By Theorem \ref{thm:Hausdorff_Young_k_transform}, we have
\begin{align*}
\|\mathcal{F}(x)\|_p\leq k^{\frac{1}{p}}\|x\|_q.
\end{align*}
For the adjoint operator $\mathcal{F}^*$, we have
\begin{align*}
\|\mathcal{F}^*\|_{1\to\infty}\leq1,\quad \|\mathcal{F}^*\|_{2\to2}=\sqrt{k}.
\end{align*}
Applying the same process in Theorem \ref{thm:Hausdorff_Young_k_transform} to $\mathcal{F}^*$, we also have
\begin{align*}
\|\mathcal{F}^*\mathcal{F}(x)\|_p\leq k^{\frac{1}{p}}\|\mathcal{F}(x)\|_q.
\end{align*}
Multiplying the above two inequalities, we obtain
\begin{align*}
\|x\|_q\|\mathcal{F}(x)\|_q\geq k^{-\frac{2}{p}}\|\mathcal{F}^*\mathcal{F}(x)\|_p\|\mathcal{F}(x)\|_p=k^{1-\frac{2}{p}}\|x\|_p\|\mathcal{F}(x)\|_p.
\end{align*}
This completes the proof of the theorem.
\end{proof}

Next, we introduce the quantum Wigderson-Wigderson uncertainty principle for quantum Fourier
transform for any $0<p,q\leq\infty$, based on the norm of quantum Fourier transform computed in \cite{LW}.
\begin{theorem}[The norm of quantum Fourier transform]\label{thm:norm_quantum_Fourier_transform}
Suppose $\mathcal{P}$ is an irreducible subfactor planar algebra. Let $\mathscr{F}$ be the Fourier transform from $\mathcal{P}_{2,\pm}$ onto $\mathcal{P}_{2,\mp}$.
Let $x\in\mathcal{P}_{2,\pm}$ be a 2-box and $0< p,q\leq \infty$. Then
\begin{align*}
K\bigg(\frac{1}{p},\frac{1}{q}\bigg)^{-1} \|x\|_q  \leq \|\mathscr{F}(x)\|_p\leq K\bigg(\frac{1}{q},\frac{1}{p}\bigg)\|x\|_q.
\end{align*}

\end{theorem}
We refer the readers to Appendix \ref{sec:appendix A} for  the specific definition of
the function  $K(\frac{1}{p},\frac{1}{q})=\|\mathcal{F}\|_{p\to q}$.

The following theorem follows immediately from  Theorem \ref{thm:norm_quantum_Fourier_transform}.
\begin{theorem}[The quantum Wigderson-Wigderson uncertainty principle for quantum Fourier transform]
Let $x\in\mathcal{P}_{2,\pm}$, we have
\begin{align*}
\|x\|_q\|\mathscr{F}(x)\|_q\geq K\bigg(\frac{1}{q},\frac{1}{p}\bigg)^{-2}\|x\|_p
\|\mathscr{F}(x)\|_p.
\end{align*}
\begin{proof}
By Theorem \ref{thm:norm_quantum_Fourier_transform}, we have
\begin{align*}
\|\mathscr{F}(x)\|_p\leq K\bigg(\frac{1}{q},\frac{1}{p}\bigg)\|x\|_q,\quad  \|x\|_p  \leq K\bigg(\frac{1}{q},\frac{1}{p}\bigg)\|\mathscr{F}(x)\|_q.
\end{align*}
Multiplying the above two equations, we can obtain the result.
\end{proof}
\end{theorem}

\subsection{Quantum Hirschman-Beckner uncertainty principle}\label{subsec:Quantum_ Hirschman_Beckner_uncertainty_principle}
In this subsection, we will prove the  quantum (smooth)
Hirschman-Beckner uncertainty principle (See Theorems \ref{thm:quantum_Hirschman_k_transform} and \ref{thm:quantum_smooth Hirschman_Beckner_k_transform}) for von Neumann
$k$-bi-algebras. For classical Hirschman-Beckner uncertainty principle \cite{Hir,Bec}, the Shannon entropy is used  to describe the uncertainty principle on $\mathbb{R}$. For finite von Neumann algebras, we would like to use von Neumann entropy instead of Shannon entropy.

\begin{definition}\label{def:von_entropy}
Let $\mathcal{M}$ be a finite von Neumann algebra with a positive  trace $\tau_{\mathcal{M}}$. The von Neumann entropy of $|x|^2\in\mathcal{M}$ is defined as follows
\begin{align*}
H(|x|^2):=-\tau_{\mathcal{M}}(|x|^2\log|x|^2)=-\tau_{\mathcal{M}}(x^*x\log x^*x).
\end{align*}
\end{definition}

We have the quantum
Hirschman-Beckner uncertainty principle for von Neumann $k$-bi-algebras.

\begin{theorem}\label{thm:quantum_Hirschman_k_transform}
Let $(\mathcal{A},\mathcal{B},d,\tau,\mathcal{F})$
be a von Neumann $k$-bi-algebra.
Suppose $\mathcal{A}$ and $\mathcal{B}$ are finite-dimensional and $\mathcal{F}^*\mathcal{F}=kI$.
 Let $x$ be a non-zero element in $\mathcal{A}$. Then we have
\begin{align*}
\frac{H(|x|^2)}{\|x\|_2^2}+\frac{H(|\mathcal{F}(x)|^2)}{\|\mathcal{F}(x)\|_2^2}\geq-\log\|x\|_2^2-\log\|\mathcal{F}(x)\|_2^2+\log k.
\end{align*}
In particular, since $\|\mathcal{F}(x)\|_2^2=k\|x\|_2^2$, we have
\begin{align*}
\frac{H(|x|^2)}{\|x\|_2^2}+\frac{H(|\mathcal{F}(x)|^2)}{\|\mathcal{F}(x)\|_2^2}\geq0
\end{align*}
whenever $\|x\|_2=1$.
\begin{proof}
By Theorem \ref{thm:Hausdorff_Young_k_transform}, we have
\begin{align*}
\|\mathcal{F}(x)\|_q\leq k^{\frac{1}{q}}\|x\|_p,
\end{align*}
where $2\leq q\leq\infty$ and $\frac{1}{p}+\frac{1}{q}=1$. Let $f(q)=\log\|\mathcal{F}(x)\|_q-\log\|x\|_p-\frac{1}{q}\log k$, then $f(q)\leq 0$ and $f(2)=0$, which implies $f'(2)\leq0$. Let $|\mathcal{F}(x)|=\sum_{i=1}^n \lambda_i e_i$ be the spectral decomposition. We have
\begin{align*}
\frac{d}{dq}\|\mathcal{F}(x)\|_q^q\bigg|_{q=2}=\frac{d}{dq}
\sum_{i=1}^n|\lambda_i|^q\tau(e_i)\bigg|_{q=2}=
\sum_{i=1}^n|\lambda_i|^2\log|\lambda_i|\tau(e_i)=-\frac{1}{2}H(|\mathcal{F}(x)|^2).
\end{align*}
Analogously,
\begin{align*}
\frac{d}{dq}\|x\|_p^p\bigg|_{q=2}=\frac{dp}{dq}\bigg|_{q=2}\frac{d}{dp}\|x\|_p^p\bigg|_{p=2}=
\frac{1}{2}H(|x|^2).
\end{align*}
Thus
\begin{align*}
\frac{d}{dq}\log\|\mathcal{F}(x)\|_q\bigg|_{q=2}=-\frac{1}{4}\log\|\mathcal{F}(x)\|_2^2-\frac{H(|\mathcal{F}(x)|^2)}{4\|\mathcal{F}(x)\|_2^2},
\end{align*}
and
\begin{align*}
\frac{d}{dq}\log\|x\|_p\bigg|_{q=2}=\frac{1}{4}\log\|x\|_2^2+\frac{H(|x|^2)}{4\|x\|_2^2}.
\end{align*}
We have
\begin{align*}
f'(2)=-\frac{1}{4}\log\|\mathcal{F}(x)\|_2^2-\frac{H(|\mathcal{F}(x)|^2)}{4\|\mathcal{F}(x)\|_2^2}-\frac{1}{4}\log\|x\|_2^2-\frac{H(|x|^2)}{4\|x\|_2^2}+\frac{1}{4}\log k.
\end{align*}
Since $f'(2)\leq0$, we obtain
\begin{align*}
\frac{H(|x|^2)}{\|x\|_2^2}+\frac{H(|\mathcal{F}(x)|^2)}{\|\mathcal{F}(x)\|_2^2}\geq-\log\|x\|_2^2-\log\|\mathcal{F}(x)\|_2^2+\log k.
\end{align*}
This completes the proof.
\end{proof}
\end{theorem}

\begin{remark}
Using the inequality $\log \mathcal{S}(x)\geq H(|x|^2)$ when $\|x\|_2=1$, we could obtain
\[\mathcal{S}(x)\mathcal{S}(\mathcal{F}(x))\geq k,\] the quantum support uncertainty principle (See Theorem \ref{thm:quantum_support-size _uncertainty_principle}).
\end{remark}
A natural question is to consider the perturbations of the inequality in Theorem  \ref{thm:quantum_Hirschman_k_transform}.
We firstly consider the smooth von Neumann entropy.
\begin{definition}\label{def:smooth_von_entropy}
Let $\mathcal{M}$ be a finite von Neumann algebra.
For any $x\in\mathcal{M}$, $\epsilon\in [0,1]$ and $p\in[1,\infty]$,
 the $(p,\epsilon)$ \textbf{smooth entropy} of $|x|^2$ is defined by
\begin{align*}
H_{\epsilon}^p(|x|^2)&:=\inf\{H(|y|^2):\ y\in\mathcal{M},\ \|x-y\|_p\leq\epsilon\},
\end{align*}
\end{definition}

\begin{remark}We thank Kaifeng Bu for referring us to
another smooth Renyi entropy studied by R. Renner and S. Wolf in quantum information in \cite{RS}.
\end{remark}
The von Neumann entropy  is continuous with respect to the operator norm and satisfies the Lipschitz condition.

\begin{lemma}\label{lem:two_matrix_spectrum_infty}
Let $A,B$ be two  matrices in $M_n(\mathbb{C})$.
Let $\lambda_1\geq\cdots\geq\lambda_n\geq0$ and $\lambda_1'\geq\cdots\geq\lambda_n'\geq0$ be the eigenvalues of
$|A|$ and $|B|$, respectively. Then we have
\[\sup_{1\leq i\leq n}|\lambda_i-\lambda_i'|\leq\|A-B\|.\]
\begin{proof}
Let $\{e_i\}_{i=1}^n$ and $\{f_i\}_{i=1}^n$ be two orthonormal basis of $\mathbb{C}^n$ such that $|A|e_i=\lambda_i e_i$ and $|B|f_i=\lambda_i' f_i$.
Since $|\|A\|-\|B\||\leq\|A-B\|$, we have $|\lambda_1-\lambda_1'|\leq\|A-B\|$. Suppose there exits $1<j\leq n$
such that $|\lambda_i-\lambda_i'|\leq\|A-B\|$ for $i<j$
and $|\lambda_j-\lambda_j'|>\|A-B\|$, we may assume that $\lambda_j>\lambda_j'+\|A-B\|$.
Let $E$ be the projection from $\mathbb{C}^n$ onto
$span\{e_1,\ldots,e_j\}$ and $F$ be the projection from $\mathbb{C}^n$ onto
$span\{f_j,\ldots,f_n\}$, then $E\wedge F\neq \emptyset$. Take a unit vector $v$ in $E\wedge F$,
we have
\begin{align*}
\|A-B\|\geq\|(A-B)v\|_2&\geq \|Av\|_2-\|Bv\|_2\\
 & = \||A|v\|_2-\||B|v\|_2\geq
\lambda_j-\lambda_j'>\|A-B\|,
\end{align*}
which leads to a contradiction.
\end{proof}

\end{lemma}

\begin{proposition}[Lipschitz condition]\label{prop:von_entropy_continuity_infty}
Let $\mathcal{M}$ be a finite-dimensional von Neumann algebra with a  trace $\tau_{\mathcal{M}}$.
For any $x,y\in\mathcal{M}$, let
$t=\max\{\|x\|,\|y\|,1\}$,
we have
\begin{align*}
|H(|x|^2)-H(|y|^2)|\leq f(t)\tau_{\mathcal{M}}(I)\|x-y\|,
\end{align*}
where $f(t)=4t\log t+2t$.
\begin{proof}
Since  $\mathcal{M}$ is finite-dimensional, we may assume  that
\[\mathcal{M}\quad=\quad   \bigoplus_{i=1}^m\mathop{M_{n_i}(\mathbb{C})}\limits_{\delta_i}.\]
Let $Tr_i$ be the unnormalized trace on $M_{n_i}(\mathbb{C})$, then we have $\tau_{\mathcal{M}}=\sum_{i=1}^m\delta_i Tr_i$. In particular, $\tau_{\mathcal{M}}(I)=\sum_{i=1}^m\delta_i n_i$.
 Let $x=\sum_{i=1}^m x_i$ and $y=\sum_{i=1}^m y_i$, let $\alpha_{i1}\geq\cdots\geq\alpha_{in_i}$
and $\beta_{i1}\geq\cdots\geq\beta_{in_i}$ be eigenvalues of $|x_i|$ and $|y_i|$, respectively.
Then we have
\begin{align*}
|H(|y|^2)-H(|x|^2)|&=|\tau_{\mathcal{M}}(|y|^2\log|y|^2-|x|^2\log|x|^2)|\\
&\leq \sum_{i=1}^m\delta_i|Tr_i(|y_i|^2\log|y_i|^2-|x_i|^2\log|x_i|^2)|\\
&=
 \sum_{i=1}^m\delta_i|\sum_{j=1}^{n_i}(\alpha_{ij}^2\log\alpha_{ij}^2-\beta_{ij}^2\log\beta_{ij}^2)|\\
 &\leq
  \sum_{i=1}^m\delta_i \sum_{j=1}^{n_i}
  |\alpha_{ij}-\beta_{ij}|\cdot(4t\log t+2t).
\end{align*}
By Lemma \ref{lem:two_matrix_spectrum_infty}, we have $|\alpha_{ij}-\beta_{ij}|\leq\|x-y\|$.  Then
\begin{align*}
|H(|y|^2)-H(|x|^2)|&\leq
  \sum_{i=1}^m\delta_i n_i\cdot
(4t\log t+2t)\|x-y\|=f(t)\tau_{\mathcal{M}}(I) \|x-y\|.
\end{align*}
\end{proof}
\end{proposition}

Suppose $\mathcal{A}$ and $\mathcal{B}$ are finite-dimensional, let
\begin{align*}
d_1&=\min\{d(e):\ e\ \text{is a minimal projection in}\ \mathcal{A} \}\\
\tau_1&=\min\{\tau(f):\ f\ \text{is a minimal projection in}\ \mathcal{B} \}.
\end{align*}
We have the quantum smooth Hirschman-Beckner uncertainty principle.
\begin{theorem}\label{thm:quantum_smooth Hirschman_Beckner_k_transform}
Let $(\mathcal{A},\mathcal{B},d,\tau,\mathcal{F})$
be a von Neumann $k$-bi-algebra.
Suppose $\mathcal{A}$ and $\mathcal{B}$ are finite-dimensional and $\mathcal{F}^*\mathcal{F}=kI$. Let $x$ be a non-zero element in $\mathcal{A}$.
 For any $\epsilon,\eta\in [0,1]$ and $p,q\in[1,\infty]$, we have
\begin{align*}
\frac{H_{\epsilon}^p(|x|^2)}{\|x\|_2^2}+\frac{H_{\eta}^q(|\mathcal{F}(x)|^2)}{\|\mathcal{F}(x)\|_2^2}\geq-4\log\|x\|_2 -\frac{C_1(x)}{\|x\|_2^2}d_1^{-\frac{1}{p}}d(I)\epsilon-
\frac{C_2(x)}{\|\mathcal{F}(x)\|_2^2}\tau_1^{-\frac{1}{q}}\tau(I)\eta,
\end{align*}
where $C_1(x)=f(\|x\|+1)$ and $C_2(x)=f(\|\mathcal{F}(x)\|+1)$
and $f(t)=4t\log t+2t$, $d_1$ and $\tau_1$ are two constants independent of $x$.
\begin{proof}
By Proposition \ref{prop:von_entropy_continuity_infty},
for any $y\in\mathcal{A}$ with $\|x-y\|_p\leq\epsilon$, we have
\begin{align*}
|H(|y|^2)-H(|x|^2)|\leq
C_1(x)d(I)\|x-y\|\leq
C_1(x)d_1^{-\frac{1}{p}}d(I)\|x-y\|_{p}.
\end{align*}
Thus we have
\[H_{\epsilon}^{p}(|x|^2)\geq
H(|x|^2)-C_1(x)d_1^{-\frac{1}{p}}d(I) \epsilon.\]
Analogously, we have
\[H_{\eta}^{q}(|\mathcal{F}(x)|^2)\geq
H(|\mathcal{F}(x)|^2)-C_2(x)\tau_1^{-\frac{1}{q}}\tau(I) \eta.
\]
Adding the above two equations and applying Theorem \ref{thm:quantum_Hirschman_k_transform}, we obtain the result.
\end{proof}
\end{theorem}

Applying Theorem \ref{thm:quantum_smooth Hirschman_Beckner_k_transform} to quantum Fourier transform on subfactor planar algebras, we have the following corollary.
\begin{corollary}\label{cor:quantum_smooth_HB_Fourier_transform}
Suppose $x$ is a non-zero $2$-box in an irreducible subfactor planar algebra. For any $\epsilon,\eta\in [0,1]$ and $p,q\in[1,\infty]$, we have
\begin{align*}
H_{\epsilon}^{p}(|x|^2)+H_{\eta}^{q}(|\mathscr{F}(x)|^2)&\geq \|x\|_2^2(2\log\delta-4\log\|x\|_2)
-   C(x)\delta^2(\epsilon+\eta),
\end{align*}
where $C(x)=f(\|x\|_2+1)$ and $f(t)=4t\log t+2t$
$4(\|x\|_2+1)\log(\|x\|_2+1)+2(\|x\|_2+1)$ and $\delta$ is the square root of Jones index.
\end{corollary}

\begin{remark}
We have the following statements:
\begin{enumerate}
\item In Example \ref{ex:fushion bi-algebra},
 take $\epsilon=\eta=0$, then Theorem \ref{thm:quantum_smooth Hirschman_Beckner_k_transform} implies Theorem 4.9 in \cite{LPW};
\item In Example \ref{ex:subfactor planar algebra},
  take $\epsilon=\eta=0$, then Theorem \ref{thm:quantum_smooth Hirschman_Beckner_k_transform} implies Theorem 5.5 in \cite{JLW}.
\end{enumerate}
\end{remark}
In \cite{JLW}, the minimizers of Hirschman-Beckner uncertainty principle on subfactor planar algebras were characterized as bi-shifts of biprojections (See Theorems 6.4 and 6.13 in \cite{JLW}). So it is natural to ask the following inverse problem.
\begin{problem}
Find a positive function $C(\epsilon,\delta)$, for $\epsilon,\delta>0$, such that $\displaystyle \lim_{\epsilon \to 0} C(\epsilon,\delta) \to 0$, and
for any 2-box $x$ of any irreducible subfactor planar algebra with Jones index $\delta^2$, $\|x\|_2=1$,
if \[H(|x|^2)+H(|\mathscr{F}(x)|^2)\geq 2\log\delta-\epsilon,\]
then $\|x-y\|\leq C(\epsilon,\delta)$ for some bi-shift of biprojection $y$.
\end{problem}

 \section{An answer to a conjecture of Wigderson and Wigderson}\label{sec:conjecture}

The famous Heisenberg uncertainty principle in \cite{Hei} could be mathematically formulated in terms of Schwarz functions on $\mathbb{R}$, (see e.g. \cite{Ken, Weyl} and Theorem 4.9 in \cite{WW}), as follows:
\begin{theorem}[Heisenberg's uncertainty principle]
Let $\mathcal{S}(\mathbb{R})$ be the space of Schwartz functions. For any $f\in\mathcal{S}(\mathbb{R})$,
\begin{align*}
\int_{\mathbb{R}}x^2|f(x)|^2dx\int_{\mathbb{R}}\xi^2|\hat{f}(\xi)|^2dx\geq\frac{1}{16\pi^2}\|f\|_2^2\|\hat{f}\|_2^2,
\end{align*}
where $\hat{f}(\xi)=\int_{\mathbb{R}}f(x)e^{-2\pi i x\xi}dx$ is the Fourier transform of $f$.
\end{theorem}
 In \cite{WW}, A. Wigderson and Y. Wigderson proved the following generalization of Heisenberg's uncertainty principle for arbitrary $q$-norm.
\begin{theorem}[See Theorem 4.11 in \cite{WW}]
 For any $f\in\mathcal{S}(\mathbb{R})$, and any $1<q\leq \infty$,
 \begin{align*}
\int_{\mathbb{R}}x^2|f(x)|^2dx\int_{\mathbb{R}}\xi^2|\hat{f}(\xi)|^2dx\geq 2^{-\frac{10q-8}{q-1}}\|f\|_q^2\|\hat{f}\|_q^2.
\end{align*}
\end{theorem}
\noindent In order to compare these inequalities for different $q$, they proposed the following conjecture
\begin{conjecture}[Conjecture 4.13 in \cite{WW}]\label{Conj: Wigderson}
For any non-zero $f \in \mathcal{S}(\mathbb{R})$, $q\in(1,\infty]$, define
\begin{align*}
F_q(f)=\frac{\|f\|_q\|\hat{f}\|_q}{\|f\|_2\|\hat{f}\|_2}=\frac{\|f\|_q\|\hat{f}\|_q}{\|f\|_2^2}.
\end{align*}
Then the image of $F_q$: $\mathcal{S}(\mathbb{R})\setminus\{0\}\to \mathbb{R}_{>0}$, is $\mathbb{R}_{>0}$ for all $q\neq 2$.
\end{conjecture}
\noindent Moreover, they proved the conjecture for $q=\infty$ in Theorem 4.12 in \cite{WW}.

In the following theorem, we verify Conjecture \ref{Conj: Wigderson} for $q>2$ and disprove Conjecture \ref{Conj: Wigderson}  for $1<q<2$. More precisely,
\begin{theorem}\label{thm:answer to conjecture}
\begin{enumerate}
\item If $1<q<2$,  take $1/p+1/q=1$, then
$$F_q(f)\geq [p^{1/p}/q^{1/q}]^{1/2}, ~\forall f\in\mathcal{S}(\mathbb{R})\setminus\{ 0\}.$$
\item If $q>2$, then the image of $F_q$ is $\mathbb{R}_{>0}$.
\end{enumerate}
\end{theorem}
To prove Theorem \ref{thm:answer to conjecture},
we firstly prove a technical lemma.
\begin{lemma}\label{lem:infty_zero_sequence}
Let $1<p<q\leq\infty$,
and define a function
$F_{p,q}$: $\mathcal{S}(\mathbb{R})\setminus\{0\}\to \mathbb{R}_{>0}$ by
\begin{align*}
F_{p,q}(f)=\frac{\|f\|_q\|\hat{f}\|_q}{\|f\|_p\|\hat{f}\|_p}.
\end{align*}
If there exist two sequences $\{f_n\}_{n=1}^{\infty}$ and $\{g_n\}_{n=1}^{\infty}$ in $\mathcal{S}(\mathbb{R})\setminus\{0\}$ such that
\begin{align*}
\lim_{n\to\infty}F_{p,q}(f_n)=0,\quad
\lim_{n\to\infty}F_{p,q}(g_n)=\infty,
\end{align*}
and
$\lambda f_n+(1-\lambda)g_n\neq0$ for any $\lambda\in[0,1]$ and all $n\geq1$,
then the image of $F_{p,q}$ is all of $\mathbb{R}_{>0}$.
 \begin{proof}
 We define
 \begin{align*}
 h_n(\lambda)=F_{p,q}(\lambda f_n+(1-\lambda)g_n),
 \end{align*}
 then $h_n(\lambda)$ is a continuous function for any $n\geq1$.
 Thus $h_n(\lambda)$ can take all real values between
 $F_{p,q}(f_n)$ and $F_{p,q}(g_n)$. The result follows immediately from the assumptions.
 \end{proof}
\end{lemma}
\begin{theorem}\label{thm:range_fourier}
Let $1<p<q\leq\infty$,
and define a function
$F_{p,q}$: $\mathcal{S}(\mathbb{R})\setminus\{0\}\to \mathbb{R}_{>0}$ by
\begin{align*}
F_{p,q}(f)=\frac{\|f\|_q\|\hat{f}\|_q}{\|f\|_p\|\hat{f}\|_p}.
\end{align*}
When $\frac{1}{p}+\frac{1}{q}<1$,
then the image of $F_{p,q}$ is all of $\mathbb{R}_{>0}$.
\begin{proof}
We will consider two special families of Schwartz functions.
Let $a>b>0$ be real numbers. Define
\begin{align*}
f_{a,b}(x)=e^{-\pi ((a+ib)x)^2}=e^{-\pi(a^2-b^2)x^2}e^{-2\pi iabx^2}.
\end{align*}
From the definition, we see that $|f_{a,b}(x)|=e^{-\pi(a^2-b^2)x^2}$. We can compute the $r$-norm of $f_{a,b}$
\begin{align*}
\quad\|f_{a,b}(x)\|_r=\bigg(\int_{-\infty}^{\infty}e^{-r\pi(a^2-b^2)x^2}dx\bigg)^{\frac{1}{r}}=
\bigg(\frac{1}{\sqrt{r(a^2-b^2)}}\bigg)^{\frac{1}{r}},\quad 1<r\leq\infty.
\end{align*}
When $r=\infty$, the above equality means $\|f_{a,b}(x)\|_{\infty}=1$.
The Fourier transform of $f_{a,b}$ is
\begin{align*}
\widehat{f_{a,b}}(\xi)=\frac{1}{a+bi}e^{-\pi(\xi/(a+bi))^2}=\frac{1}{a+bi}e^{-\pi\xi^2(a^2-b^2)/(a^2+b^2)^2}e^{2\pi i\xi^2ab/(a^2+b^2)^2}.
\end{align*}
In particular,
\begin{align*}
\big|\widehat{f_{a,b}}(\xi)\big|=\frac{1}{\sqrt{a^2+b^2}}e^{-\pi\xi^2(a^2-b^2)/(a^2+b^2)^2}.
\end{align*}
Similarly, we can compute the $r$-norm of $\widehat{f_{a,b}}(\xi)$
\begin{align*}
\|\widehat{f_{a,b}}(\xi)\|_r=\frac{1}{\sqrt{a^2+b^2}}\bigg(\frac{a^2+b^2}{\sqrt{r(a^2-b^2)}}\bigg)^{\frac{1}{r}},\quad 1<r\leq\infty.
\end{align*}
When $r=\infty$, $\|\widehat{f_{a,b}}(\xi)\|_{\infty}=\frac{1}{\sqrt{a^2+b^2}}$.
This implies that
\begin{align*}
F_{p,q}(f_{a,b})=\frac{\|f\|_q\|\hat{f}\|_q}{\|f\|_p\|\hat{f}\|_p}=\frac{\sqrt[p]{p}(a^2+b^2)^{\frac{1}{q}-\frac{1}{p}}}{ \sqrt[q]{q} (a^2-b^2)^{\frac{1}{q}-\frac{1}{p}}},\quad 1<p<q\leq\infty.
\end{align*}
If $a>1$ and $b=\sqrt{a^2-1}$, then we get
\begin{align*}
F_{p,q}(f_{a,\sqrt{a^2-1}})=\frac{\sqrt[p]{p}}{ \sqrt[q]{q} }(2a^2-1)^{\frac{1}{q}-\frac{1}{p}}.
\end{align*}
So
\begin{align}\label{eq:f_ab_infty}
\lim_{a\to\infty}F_{p,q}(f_{a,\sqrt{a^2-1}})=
0,\quad 1<p<q\leq\infty.
\end{align}

Next, we consider another family of functions
\begin{align*}
g_c(x)=\frac{1}{\sqrt{c}}e^{-\pi(x/c)^2}+\sqrt{c}e^{-\pi(cx)^2}, \quad c>0 \;.
\end{align*}
Both items of $g_c$ are positive, so
\begin{align*}
\frac{1}{2}(\frac{1}{\sqrt{c} }\|e^{-\pi(x/c)^2} \|_r+
\sqrt{c}\|e^{-\pi(cx)^2}\|_r)<\|g_c\|_r<
\frac{1}{\sqrt{c} }\|e^{-\pi(x/c)^2} \|_r+
\sqrt{c}\|e^{-\pi(cx)^2}\|_r
\end{align*}
for any $1<r\leq\infty$, namely
\begin{align*}
\frac{c^{\frac{1}{r}-\frac{1}{2}}+c^{\frac{1}{2}-\frac{1}{r}}}{2\sqrt[2r]{r}}<\|g_c\|_r<\frac{c^{\frac{1}{r}-\frac{1}{2}}+c^{\frac{1}{2}-\frac{1}{r}}}{\sqrt[2r]{r}} \;.
\end{align*}
Note that $\hat{g_c}=g_c$, so we obtain an estimation of $F_{p,q}(g_c)$
\begin{align*}
\frac{\sqrt[p]{p}(c^{\frac{1}{q}-\frac{1}{2}}+c^{\frac{1}{2}-\frac{1}{q}})^2}{4\sqrt[q]{q}(c^{\frac{1}{p}-\frac{1}{2}}+c^{\frac{1}{2}-\frac{1}{p}})^2}<F_{p,q}(g_c)<\frac{4\sqrt[p]{p}(c^{\frac{1}{q}-\frac{1}{2}}+c^{\frac{1}{2}-\frac{1}{q}})^2}{\sqrt[q]{q}(c^{\frac{1}{p}-\frac{1}{2}}+c^{\frac{1}{2}-\frac{1}{p}})^2},\quad 1<p<q\leq\infty.
\end{align*}
Therefore,
\begin{align}\label{eq:g_c_infty}
\lim_{c\to\infty}F_{p,q}(g_c)=
\begin{cases}
\infty,&\frac{1}{p}+\frac{1}{q}<1,\ 1<p<q\leq\infty,\\
0,&\frac{1}{p}+\frac{1}{q}>1,\ 1<p<q\leq\infty.
\end{cases}
\end{align}
Combining equation (\ref{eq:f_ab_infty}) and equation (\ref{eq:g_c_infty}), we have
\begin{align*}
\lim_{a\to\infty}F_{p,q}(f_{a,\sqrt{a^2-1}})=0,\quad\lim_{c\to\infty}F_{p,q}(g_c)=\infty
\end{align*}
when $\frac{1}{p}+\frac{1}{q}<1$.
Then by Lemma \ref{lem:infty_zero_sequence}, we can obtain the conclusion.
\end{proof}
\end{theorem}
In particular, take $p=2$, then $2<q\leq\infty$, which implies that Conjecture \ref{Conj: Wigderson} holds for all $2<q\leq\infty$.
\begin{theorem}\label{thm:p=2,1<q<2,lower bound}
When $1<q<2$, we have
\begin{align*}
F_q(f)=\frac{\|f\|_q\|\hat{f}\|_q}{\|f\|_2\|\hat{f}\|_2}=\frac{\|f\|_q\|\hat{f}\|_q}{\|\hat{f}\|_2^2}\geq [p^{1/p}/q^{1/q}]^{1/2},
\end{align*}
where $ \frac{1}{p}+\frac{1}{q}=1$,
for any $f\in\mathcal{S}(\mathbb{R})\setminus\{0\}$.
\begin{proof}
By the sharp Hausdorff-Young inequality (See Theorem 1 in
\cite{Bec}), we have
\begin{align*}
\|\hat{f}\|_p\leq [q^{1/q}/p^{1/p}]^{1/2} \|f\|_q
\end{align*}
for any $f\in\mathcal{S}(\mathbb{R})\setminus\{0\}$.
By H\"older's inequality, we further have
\begin{align*}
\|\hat{f}\|_2^2\leq\|\hat{f}\|_p\|\hat{f}\|_q.
\end{align*}
Combining the above two equations, we obtain
\begin{align*}
F_q(f)\geq [p^{1/p}/q^{1/q}]^{1/2}\frac{\|\hat{f}\|_p\|\hat{f}\|_q}{\|\hat{f}\|_2^2}\geq [p^{1/p}/q^{1/q}]^{1/2} .
\end{align*}
\end{proof}
\end{theorem}
Now, we give the proof of Theorem \ref{thm:answer to conjecture}.
\begin{proof}[Proof of Theorem \ref{thm:answer to conjecture}]
Combining Theorem \ref{thm:range_fourier} and Theorem \ref{thm:p=2,1<q<2,lower bound}, we can obtain the conclusion.
\end{proof}
By Theorem \ref{thm:answer to conjecture}, we have known that $F_q(f)\geq [p^{1/p}/q^{1/q}]^{1/2}$ for all $f\in\mathcal{S}(\mathbb{R})\setminus\{ 0\}$ when $1<q<2$. Let $C_q=\inf\{F_q(f):\ f\in\mathcal{S}(\mathbb{R})\setminus\{0\} \}$, then $C_q\geq [p^{1/p}/q^{1/q}]^{1/2}$ by Theorem \ref{thm:p=2,1<q<2,lower bound}. So it is natural to ask what the optimal constant $C_q$ is.
\begin{problem}
Determine the constant $C_q=\inf\{F_q(f):\ f\in\mathcal{S}(\mathbb{R})\setminus\{0\} \}$ when $1<q<2$.
\end{problem}

%

\begin{appendix}
    \renewcommand{\thesection}{\Alph{section}}
      \section{The function K(1/p,1/q)}\label{sec:appendix A}

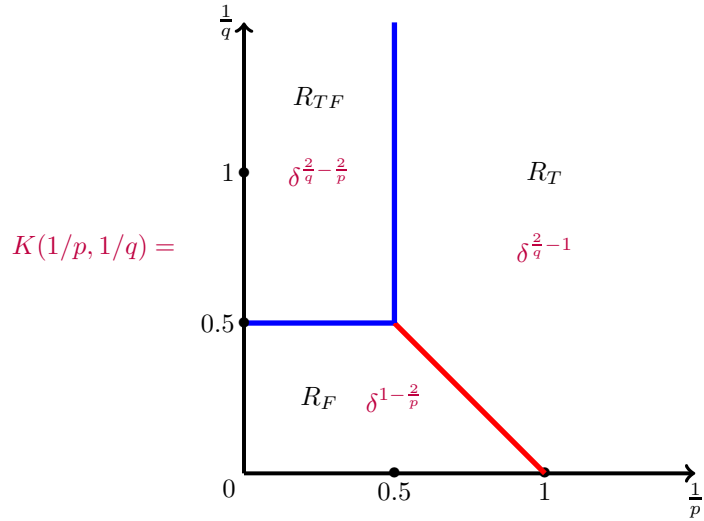
\begin{figure}[h]
\begin{center}
\begin{tikzpicture}[scale=2]
\node at (-1,1.5) {$\textcolor{purple}{K(1/p,1/q)=}$};
\draw [line width=1.5] [->](0,0)--(0,3);
\node at (-0.1,-0.1) {0};
\node [left] at (0, 3) {$\frac{1}{q}$};
\draw [line width=1.5]  [->](0,0)--(3,0);
\node [below] at (1,0) {$0.5$};
\node at (1,0) {$\bullet$};
\node [below] at (2, 0) {1};
\node at (2,0) {$\bullet$};
\node [below] at (3,0) {$\frac{1}{p}$};
\draw [blue, line width=2] (1,1)--(0,1);
\node [left] at (0,1) {$0.5$};
\node at (0,1) {$\bullet$};
\node [left] at (0, 2) {1};
\node at (0,2) {$\bullet$};
\node [purple] at (0.5,2) {$\delta^{\frac{2}{q}-\frac{2}{p}}$};
\node at (0.5,2.5) {$R_{TF}$};
\draw [red, line width=2] (1,1)--(2,0);
\node at (0.5, 0.5) {$R_F$};
\node [purple] at (1, 0.5) {$\delta^{1-\frac{2}{p}}$};
\draw [blue, line width=2] (1,1)--(1,3);
\node [purple] at (2, 1.5) {$\delta^{\frac{2}{q}-1}$};
\node at (2,2) {$R_T$};
\end{tikzpicture}
\end{center}
\caption{The norm of the Fourier transform $\mathscr{F}$.}
\label{Fig:SFT norm}
\end{figure}
The first quadrant is divided into three regions $R_{T},R_{F},R_{TF}$ as follows:
\begin{eqnarray*}
R_F:&=&\{(1/p,1/q)\in [0,\infty]^2: 1/p+1/q\leq 1,1/q\leq 1/2\},\\
R_T:&=&\{(1/p,1/q)\in [0,\infty]^2: 1/p+1/q\geq  1,1/p\geq 1/2\},\\
R_{TF}:&=&\{(1/p,1/q)\in [0,\infty]^2: 1/p\leq 1/2, 1/q\geq 1/2\}.
\end{eqnarray*}
The function $K(1/p.1/q)$ on $[0,\infty)^2$ is given by
\begin{equation*}\label{K}
K(1/p,1/q)=\left\{
\begin{array}{lll}
\delta^{1-2/p} &\text{ for } &(1/p,1/q)\in R_F,\\
\delta^{2/q-1} &\text{ for } &(1/p,1/q)\in R_T,\\
\delta^{2/q-2/p} &\text{ for }& (1/p,1/q)\in R_{TF}.
\end{array}
\right.
\end{equation*}

    \end{appendix}
 \bibliographystyle{plain}

\end{document}